\documentclass[11pt]{article}
\usepackage{graphicx}
\usepackage{latexsym}
\usepackage[all]{xy}
\usepackage{amsfonts}
\usepackage{amsthm}
\usepackage{amsmath}
\usepackage{amssymb}

\def\<{{\langle}}
\def\>{{\rangle}}

\def\note#1{{}}

\def\note#1{}

\def\rend#1#2{{{\rm End}\sb{#1}(#2)}}

\def\beq{\begin{equation}}
\def\eeq{\end{equation}}

\def\id{{I}}

\def\ot{{\otimes}}

\newcounter{zlist}

\newtheorem{thm}{Theorem}[section]

\def\Label#1{\label{#1}\ifmmode\llap{[#1] }\else 
\marginpar{\smash{\hbox{\tiny [#1]}}}\fi}
\def\Label{\label}

\newtheorem{proposition}{Proposition}[section]

\newtheorem{theorem}[proposition]{Theorem}

\theoremstyle{definition}
\newtheorem{definition}[proposition]{Definition}

\theoremstyle{remark}
\newtheorem{remark}[proposition]{Remark}

\newcounter{c}

\newcommand{\etyk}[1]{\vspace{-7.4mm}$$\begin{equation}\Label{#1}
\addtocounter{c}{1}}
\renewcommand{\]}{\ifnum \value{c}=1 $$\else \end{equation}\fi}
\usepackage{fancyhdr}
\begin{document}

\vspace*{10mm}

\begin{center}
\uppercase{\textbf{YANG--BAXTER OPERATORS FROM ALGEBRA STRUCTURES AND 
LIE SUPER-ALGEBRA STRUCTURES}}
\end{center}

\bigskip

\begin{center}
\textsc{F. F. Nichita, and Bogdan Popovici}
\end{center}

\bigskip

\textsc{Abstract.}

The concept of symmetry plays an important role in solving equations and 
systems of equations.
We will present solutions for the (constant and spectral-parameter) 
Yang-Baxter equations and
Yang-Baxter systems arising from algebra structures and discuss about
their symmetries.
In the last section,
we present enhanced
versions of Theorem 1 (from
\cite{tbh}), solutions for the classical Yang-Baxter equation, and 
solutions for the Yang-Baxter systems
from Lie (super)algebras.

\bigskip


$$
\textsc{1. Introduction and preliminaries }
$$

The concept of symmetry appears in almost every scientific and artistic area.
It plays an important role in solving equations and systems of 
equations. Sometimes an equation has a high degree of symmetry, but it is
not perfectly symmetric. Similar situations, appearing in a geometric 
framework,
were considered by \cite{ag}, where concepts like ``distance from symmetry
in shape'', ``closest symmetric shape'' and ``symmetry distance'' were 
introduced.

In this paper we reveal abstract situations in which equations and some of their solutions are symmetric or almost symmetric. Besides these forms of symmetry,
we encounter the super-symmetry. The super-symmetry property is a concept
of interest at CERN, and it states that any particle has an associated super-symmetric particle (also called its super-partner). While this property has not been confirmed yet, there exist plenty of examples and applications of 
super-symmetric structures in particle physics and quantum groups.

The quantum Yang--Baxter equation (QYBE) 
first appeared in theoretical physics (\cite{yan}) and 
statistical mechanics (\cite{rba, bax}).
It 
plays a crucial role  in analysis
of integrable systems, in quantum and statistical mechanics,
in knot theory, and also in 
the
theory of quantum groups. 
On the other hand, the 
theory of integrable Hamiltonian
systems makes great use of the solutions of the one-parameter form of the
Yang-Baxter equation, since coefficients of the power series expansion of
such a solution give rise to commuting integrals of motion.
The two-parameter form of the QYBE is related to Yang's paper (\cite{yan}), and
its solutions are referred to as a colored {\em Yang-Baxter operator}.
{\em Yang--Baxter systems} emerged from the study of quantum
integrable systems, as generalizations of the QYBE related to 
nonultralocal models.

This paper presents results
on Yang-Baxter operators from algebra structures
and related topics (colored Yang-Baxter operators, Yang-Baxter
systems, Yang-Baxter operator from Lie super-algebras).
In the last section,
we present  enhanced
versions of Theorem 1 (from
\cite{tbh}), solutions for the classical Yang-Baxter equation (CYBE),
and solutions for a $WXZ$-system.

The Yang-Baxter equations 
(QYBE, CYBE, set-theoretical Yang-Baxter equation, etc) 
have some kind of symmetries, which can be used
to find solutions for them; many times, it is possible to obtain larger
classes of solutions by some sort of deformation (quantization) of those solutions.
The Yang-Baxter equations are
related to symmetric spaces (\cite{bw}), Boolean algebras (\cite{n2}),
Jordan algebras (\cite{i}), etc.
The following is a short
bibliography on
QYBE (\cite{LamRad:Int},
\cite{Kas:Qua},
\cite{GonVes:yan},
\cite{tbh}, \cite{n1}), 
 and Yang-Baxter 
systems (\cite{HlaSno:sol},
\cite{HlaKun:qua},
\cite{bn},
\cite{np},
\cite{Vla:met}).

\bigskip

\addtocounter{section}{1}

Throughout this paper $ k $ is a field. 

All tensor products appearing in this paper are defined over $k$.

For $ V $ a $ k$-space, we denote by
$ \   \tau : V \otimes V \rightarrow V \ot V \  $ the twist map defined by $ \tau (v \ot w) = w \ot v $, and by $ I: V \rightarrow V $
the identity map of the space V.

We use the following notations concerning the Yang-Baxter equation.

If $ \  R: V \ot V \rightarrow V \ot V  $
is a $ k$-linear map, then
$ {R^{12}}= R \ot I , {R^{23}}= I \ot R ,
{R^{13}}=(I\ot\tau )(R\ot I)(I\ot \tau ) $.

\begin{definition}
  An invertible  $ k$-linear map  $ R : V \ot V \rightarrow V \ot V $
  is called a Yang-Baxter
  operator if it satisfies the  equation
  \begin{equation}  \label{ybeq}
    R^{12}  \circ  R^{23}  \circ  R^{12} = R^{23}  \circ  R^{12}  \circ  R^{23}
  \end{equation}
\end{definition}

\begin{remark}
  The equation (\ref{ybeq}) is usually called the braid equation. 
  It has some kind of symmetry.
\end{remark}

\begin{remark}
  The operator $R$ satisfies (\ref{ybeq}) if and only if $R\circ \tau  $ 
  satisfies   the constant QYBE (if and only if $ \tau \circ R $ satisfies
  the constant QYBE):
  \begin{equation}   \label{ybeq2}
    R^{12}  \circ  R^{13}  \circ  R^{23} = R^{23}  \circ  R^{13}  \circ  R^{12}
  \end{equation}

  This equation is symmetric with regard to ``=''.
\end{remark}

\begin{remark}
  (i) $ \   \tau : V \otimes V \rightarrow V \ot V \  $ is an example of a 
  Yang-Baxter operator.
  
  (ii) An exhaustive list of invertible solutions for (\ref{ybeq2}) in dimension 
  2 is given in \cite{hi} and in the appendix of \cite{HlaSno:sol}.
  
  (iii) Finding all Yang-Baxter operators in dimension greater than 2 is an 
  unsolved problem.
\end{remark}

\bigskip

Let $A$ be a (unitary) associative $k$-algebra, and $ \alpha, \beta, \gamma \in k$. 
We define the
$k$-linear map:
$ \  R^{A}_{\alpha, \beta, \gamma}: A \ot A \rightarrow A \ot A, \ \ 
R^{A}_{\alpha, \beta, \gamma}( a \ot b) = \alpha ab \ot 1 + \beta 1 \ot ab -
\gamma a \ot b $.

\begin{thm}(S. D\u asc\u alescu and F. F. Nichita, \cite{DasNic:yan})
  Let $A$ be an associative 
  $k$-algebra with $ \dim A \ge 2$, and $ \alpha, \beta, \gamma \in k$. Then $ R^{A}_{\alpha, \beta, \gamma}$ is a Yang-Baxter operator if and only if one
  of the following holds:
  
  (i) $ \alpha = \gamma \ne 0, \ \ \beta \ne 0 $;
  
  (ii) $ \beta = \gamma \ne 0, \ \ \alpha \ne 0 $;
  
  (iii) $ \alpha = \beta = 0, \ \ \gamma \ne 0 $.
  
  If so, we have $ ( R^{A}_{\alpha, \beta, \gamma})^{-1} = 
  R^{A}_{\frac{1}{ \beta}, \frac{1}{\alpha}, \frac{1}{\gamma}} $ in cases (i) and
  (ii), and $ ( R^{A}_{0, 0, \gamma})^{-1} = 
  R^{A}_{0, 0, \frac{1}{\gamma}} $ in case (iii).
\end{thm}

\begin{remark} The Yang--Baxter equation plays an important role in knot
 theory. Turaev has described a general scheme to derive an invariant of 
 oriented links from a Yang--Baxter operator, provided this 
 one can be ``enhanced''.
 In \cite{mn}, we considered the problem of applying Turaev's method to the
 Yang--Baxter operators derived from algebra structures presented in
 the above theorem. We concluded that Turaev's procedure invariably 
 produces from any of those enhancements the Alexander polynomial of knots. 
\end{remark}

\bigskip

We now present the matrix form of the operator obtained in the case (i)
of the previous theorem,
$R = R^{A}_{\alpha, \beta, \alpha}: A \ot A \rightarrow A \ot A, \ \ 
R( a \ot b) = \alpha ab \ot 1 + \beta 1 \ot ab -
\alpha a \ot b $.
We consider the algebra 
$ A = \frac{ k[X]}{(X^2-mX-n)} $, where $ m,n $ are scalars.
Then $A$ has the basis $\{1, x \}$, where $x$ is the image of $X$ in the factor
ring. 

In matrix form, this operator reads:
\begin{equation} \label{rmatcon}
  \begin{pmatrix}
    \beta & 0 & 0 & 0\\
    0 & \beta - \alpha & \alpha & 0\\
    0 & 0  & \beta & 0\\
    ( \alpha + \beta ) n & \beta m & \alpha m & - \alpha
  \end{pmatrix}
\end{equation}
Let us observe that $ R'= R \circ \tau $ is a solution for the equation
(\ref{ybeq2}).
It is convenient to get rid of the auxiliary parameters and to consider
the simplest form of $ R' $:
\begin{equation} \label{rmatcon2}
  \begin{pmatrix}
    1 & 0 & 0 & 0\\
    0 & 1 & 0 & 0\\
    0 & 1-q  & q & 0\\
    \eta & 0 & 0 & -q
  \end{pmatrix}
\end{equation}
where $ \eta \in \{ 0, \ 1 \} $, and $q \in k - \{ 0 \}$.
The matrix form (\ref{rmatcon2}) was obtained as a consequence of the fact
that isomorphic algebras produce isomorphic Yang-Baxter operators.

If $q =1 , \ \  \eta =0$, the above matrix is symmetric with regard to the first 
diagonal, and can be obtained from the self-inverse operators described in
\cite{n1}.
Thus, the general case could be seen as a deformation of that case.

$$
\textsc{2. The two-parameter form of the QYBE }
$$
\addtocounter{section}{1}

Formally, a colored Yang-Baxter operator is defined as a function $ R
:X\times X \to \rend k {V\otimes V}, $ where $X$ is a set and $V$ is a
finite dimensional vector space over a field $k$. 
Thus, for any $u,v\in X$,
$R(u,v) : V\otimes V\to V\otimes V$ is a linear operator. 
We consider three operators acting on a triple
tensor product $V\otimes V\otimes V$, $R^{12}(u,v) = R(u,v)\otimes
\id$, $R^{23}(v,w)= \id\otimes R(v,w)$, and similarly $R^{13}(u,w)$ as
an operator that acts non-trivially on the first and third factor in
$V\otimes V\otimes V$. 

If $R$
satisfies the two-parameter form of the QYBE:
\begin{equation}\label{yb} 
  R^{12}(u,v)R^{13}(u,w)R^{23}(v,w) = R^{23}(v,w)
  R^{13}(u,w)R^{12}(u,v)
\end{equation} 
$ \forall \ u,v,w\in X$, then it is called 
a { \em colored Yang-Baxter operator}.

\bigskip

\begin{theorem} (F. F. Nichita and D. Parashar, \cite{np})
  Let $A$ be an associative 
  $k$-algebra with $ \dim A \ge 2$, and
  $ X \subset k $. Then,
  for any two parameters $p,q\in k$, the function
  $R:X\times X\to \rend k {A\otimes A}$ defined by
  \begin{equation}\label{rsol} 
    R(u,v)(a\otimes b) =R_{p,q}(u,v)(a\otimes b)=
    p(u-v)1\otimes ab + q(u-v)ab\otimes 1 -(pu-qv)b\otimes a,
  \end{equation}
  satisfies the colored QYBE (\ref{yb}).
\end{theorem}

\begin{remark}
  \begin{enumerate}
  \item The solution (\ref{rsol}) is related to Yang's paper  \cite{yan}.
    
  \item If $ \ pu \neq qv $ and $ \ qu \neq pv $ then the operator
    ($\ref{rsol}$) is invertible. Moreover, the following formula holds:

    $R^{-1}(u,v)(a\otimes b) = \frac{p(u-v)}{(qu-pv)(pu-qv)}ba\otimes 1 + 
    \frac{q(u-v)}{(qu-pv)(pu-qv)}1\otimes ba - \frac{1}{(pu-qv)}b\otimes a $.
  \item It follows an almost symmetric relation:

    $  \ \ R_{q,p}(u,v) = (qu-pv)(pu-qv) R_{p,q}^{-1}(u,v) \circ \tau $.
  \item     Let us consider the Theorem 2.1. 
    If we let v=0 and u=1, we obtain the 
    operator
    $ \ R(a\otimes b) =p 1\otimes ab + q ab\otimes 1 - p b\otimes a$,
    which satisfies the constant QYBE (\ref{ybeq2}).
    Notice that $ \tau \circ R $ is the Yang-Baxter operator from the 
    Theorem 1.4, case (i).
    
  \item The system of equations (\ref{e1}--\ref{e5}) is related to the 
    above theorem. It is an open problem to classify its solutions. This 
    system of equations has some remarkable symmetry properties which can 
    be used to find some solutions for it. For example, the equations 
    (\ref{e2}) and (\ref{e5}) are in some sense dual to each other.
    Likewise, (\ref{e3}) and (\ref{e4}) are in some sense dual to each 
    other.
    \begin{eqnarray} &&
      (\beta(v,w)-\gamma(v,w))(\alpha(u,v)\beta(u,w) -
      \alpha(u,w)\beta(u,v) +
      (\alpha(u,v)-\gamma(u,v)))\nonumber \\ &&\quad \quad \quad
      (\alpha(v,w)\beta(u,w) - \alpha(u,w)\beta(v,w))
      = 0 \label{e1} \\ \nonumber \\ &&
      \beta(v,w)(\beta(u,v)-\gamma(u,v))(\alpha(u,w)-\gamma(u,w)) \nonumber \\
      &&\quad \quad \quad +
      (\alpha(v,w)-\gamma(v,w))(\beta(u,w)\gamma(u,v)-\beta(u,v)\gamma(u,w)) =
      0 \label{e2} \\ \nonumber \\ && \alpha(u,v)
      \beta(v,w)(\alpha(u,w)-\gamma(u,w)) + \alpha(v,w)\gamma(u,w)
      (\gamma(u,v) - \alpha(u,v)) \nonumber \\ &&\quad \quad \quad +
      \gamma(v,w) (\alpha(u,v)\gamma(u,w)-\alpha(u,w)\gamma(u,v)) = 0
      \label{e3} \\ \nonumber \\ && \alpha(u,v)
      \beta(v,w)(\beta(u,w)-\gamma(u,w)) + \beta(v,w)\gamma(u,w) (\gamma(u,v)
      - \beta(u,v)) \nonumber \\ &&\quad \quad \quad + \gamma(v,w)
      (\beta(u,v)\gamma(u,w)-\beta(u,w)\gamma(u,v)) = 0 \label{e4} \\
      \nonumber \\ && \alpha(u,v)(\alpha(v,w)-\gamma(v,w))(\beta(u,w) -
      \gamma(u,w))  + 
      (\beta(u,v)-
      \gamma(u,v))\nonumber \\ &&\quad \quad \quad
      ( \alpha(u,w) \gamma(v,w) - \alpha(v,w) \gamma(u,w)) = 0 \
      \label{e5} 
    \end{eqnarray}
  \end{enumerate}
\end{remark}

\bigskip

We now consider the algebra 
$ A = \frac{ k[X]}{(X^2-\sigma)} $, where $ \sigma $ is a scalar.
Then $A$ has the basis $\{1, x \}$, where $x$ is the image of $X$ in the factor
ring. We consider the basis $ \{ 1 \otimes 1, 1 \otimes x, x \otimes 1, 
x \otimes x \} $ of $ A \otimes A $ and represent the operator (\ref{rsol})
in this basis:

$ R(u,v)(1\otimes 1) = (qu-pv) 1 \otimes 1 $

$ R(u,v)(1\otimes x) = p(u-v) 1 \otimes x + (q-p)u x \otimes 1 $

$ R(u,v)(x\otimes 1) = (q-p)v 1 \otimes x + q(u-v) x \otimes 1 $

$ R(u,v)(x\otimes x) = \sigma (p+q)(u-v) 1 \otimes 1 - 
(pu-qv) x \otimes x $

In matrix form, this operator reads
\begin{equation} \label{rmat}
  R(u,v)= \begin{pmatrix}
    qu-pv & 0 & 0 & \sigma (q+p)(u-v)\\
    0 & p(u-v) & (q-p)v & 0\\
    0 & (q-p)u & q(u-v) & 0\\
    0 & 0 & 0 & qv-pu
  \end{pmatrix}
\end{equation}

$$
\textsc{3. Yang-Baxter systems }
$$
\addtocounter{section}{1}

It is convenient to describe the Yang-Baxter systems in terms of
the Yang-Baxter commutators.

Let $V$, $V'$, $V''$ be finite dimensional
vector spaces over the field $k$, and let $R: V\ot
V' \rightarrow V\ot V'$, $S: V\ot V'' \rightarrow V\ot V''$ and $T:
V'\ot V'' \rightarrow V'\ot V''$ be three linear maps.
The {\em Yang--Baxter
commutator} is a map $[R,S,T]: V\ot V'\ot V'' \rightarrow V\ot V'\ot
V''$ defined by \beq [R,S,T]:= R^{12} S^{13} T^{23} - T^{23} S^{13}
R^{12}. \eeq 
Note that $[R,R,R] = 0$ is just a short-hand
notation for writing the constant QYBE (\ref{ybeq2}).

A system of linear
maps
$W: V\ot V\ \rightarrow V\ot V,\quad Z: V'\ot V'\ \rightarrow V'\ot
V',\quad X: V\ot V'\ \rightarrow V\ot V',$ is called a
$WXZ$--system if the
following conditions hold: \beq \label{ybsdoub} [W,W,W] = 0 \qquad
[Z,Z,Z] = 0 \qquad [W,X,X] = 0 \qquad [X,X,Z] = 0\eeq 
It
was observed that $WXZ$--systems with invertible $W,X$ and $Z$ can
be used to construct dually paired bialgebras of the FRT type
leading to quantum doubles. The above is one type of a constant
Yang--Baxter system that has recently been studied in \cite{np} and
also shown to be closely related to entwining structures \cite{bn}.
Thus, the $WXZ$--systems, because of their symmetry, are 
related to  ``gluing procedures'' (obtaining a bigger object from two
objects of the same kind). In the next theorem if we let $ \ \mu, \ \lambda \rightarrow 1$,
we obtain the trivial solution of a $WXZ$--system.

\bigskip

\begin{theorem} (F. F. Nichita and D. Parashar, \cite{np})
Let $A$ be a $k$-algebra, and $ \lambda, \mu \in k$. The following is a 
$WXZ$--system:

$ W : A \ot A \rightarrow A \ot A, \ \ 
W(a \ot b)= ab \ot 1 + \lambda 1 \ot ab - b \ot a $,

$ Z : A \ot A \rightarrow A \ot A, \ \ 
Z(a \ot b)= \mu ab \ot 1 +  1 \ot ab - b \ot a $,

$ X : A \ot A \rightarrow A \ot A, \ \ 
X(a \ot b)= ab \ot 1 +  1 \ot ab - b \ot a $.

\end{theorem}

\newpage
$$
\textsc{4. Applications and conclusions }
$$
\addtocounter{section}{1}

$$
\textsc{4.1 Improving theorems }
$$

Using the techniques from above we now present enhanced
versions of Theorem 1 (from
\cite{tbh}).\\

\begin{theorem} (F. F. Nichita and B. P. Popovici, \cite{nipo})
Let $V = W \oplus kc $ be a $k$-space, 
and $ f, g : V \ot V \rightarrow V $ $k$-linear maps such that
$ f, g = 0 $ on $ V \ot c + c \ot V $.
Then,
$ R: V \ot V \rightarrow V \ot V, \ 
R(v \ot w)= f(v \ot w) \ot c + c \ot g(v \ot w) $
is a solution for QYBE (\ref{ybeq2}).

\end{theorem}

\bigskip

$$
\textsc{4.2 Lie super-algebras }
$$

In particle physics, super-symmetry is a symmetry that relates elementary
particles of one spin to other particles that differ by half a unit of spin
and are known as super-partners.

According to the spin-statistic theorem, bosonic fields commute while fermionic fields anti-commute. Combining the two kinds of fields into a single algebra
requires the introduction of a  
 $ \mathbb{Z}_2$-grading under which
the bosons are the even elements and the fermions are the odd elements.
Such an algebra is called a Lie super-algebra.

\bigskip

Let $ ( L , [,] )$ be a Lie super-algebra over $k$,
and  $ Z(L) = \{ z \in L : [z,x]=0 \ \ \forall \ x \in L \} $.

For $ z \in Z(L), \ \vert z \vert =0 $ and $ \alpha \in k $ we define:

$$ { \phi }^L_{ \alpha} \ : \ L \ot L \ \ \longrightarrow \ \  L \ot L $$

$$ 
x \ot y \mapsto \alpha [x,y] \ot z + (-1)^{ \vert x \vert \vert y \vert } y \ot x \ . $$

Its inverse is:

$$ {{ \phi }^L_{ \alpha}}^{-1} \ : \ L \ot L \ \ 
\longrightarrow \ \  L \ot L $$

$$ 
x \ot y \mapsto z \ot [x,y] \ot z + \frac{1}{\alpha} (-1)^{ \vert x \vert \vert y \vert } y \ot x \ . $$

\begin{theorem} (F. F. Nichita and B. P. Popovici, \cite{nipo2})

Let  $ ( L , [,] )$ be a Lie superalgebra 
and
$ z \in Z(L), \vert z \vert = 0  $, and $ \alpha \in k $. Then:
$ \ \ \ \  { \phi }^L_{ \alpha} $ is a YB operator.
\end{theorem}

$$
\textsc{4.3 CYBE}
$$
\addtocounter{section}{1}

\begin{theorem} Let  $ ( L , [,] )$ be a Lie algebra 
and
$ z \in Z(L)$. Then:

$ r: L \ot L \ \ \longrightarrow \ \  L \ot L, \ \  
x \ot y \mapsto [x,y] \ot z - \alpha x \otimes y $ 

satisfies the classical Yang-Baxter equation:

$ [ r^{12},\  r^{13} ] \  + \  [r^{12}, \  r^{23}] \ + \  [r^{13}, \ r^{23}] = 0 $.
\end{theorem} 

{ \bf Proof.} It is left for the reader.

$$
\textsc{4.4 Poisson superalgebra}
$$

In mathematics, a Poisson superalgebra is a $Z_2$-graded 
generalization of a Poisson algebra. Specifically, a Poisson superalgebra is an (associative) superalgebra $A$ with a 
Lie superbracket  $ \{ ,\} : A \otimes A \rightarrow A $, such that $(A, \{ , \} )$ is a Lie superalgebra and the operator
$ \{ x ,\} : A \rightarrow A $ is a superderivation of A:

   $ \{ x,yz \} = \{ x,y \} z + (-1)^{|x||y|} y \{ x,z \}.$

This is one possible way of "super"izing the Poisson algebra. This gives the classical dynamics of fermion fields and classical spin-1/2 particles. The other is to define an antibracket algebra instead. This is used in the BRST and Batalin-Vilkovisky formalism.

\begin{thm}

Let $A$  be a Poisson superalgebra with a unity,  $1=1_A$, for the product *,  
such that
$ \{ x,\ 1_A \} = 0 \ \  \forall x \in A$.
Then, we have the following $WXZ$-system:

$W(x\otimes y) = \{ x,\ y \} \ot 1 \  + \  (-1)^{|x||y|}  x \ot y$;

$X(x\otimes y) = 1 \ot  \{ x,\ y \} \  + \   (-1)^{|x||y|}  x \ot y$;

$Z(x\otimes y) = 1 \ot  x*y \  + \  x*y \ot 1 \  - \  y \ot  x$.

\end{thm}

{ \bf Proof.} The condition $[W, W, W] = 0$ follows from Theorem 4.2, and $[Z, Z, Z] = 0$ 
follows from Theorem 1.5.

The condition $[X, X, Z] = 0$ is satisfied because $\{x,\}: A \to A$ is a superderivation of A.
$[W, X, X] = 0$ is a consequence of the Lie superalgebra axioms.

\bigskip


\begin{center}
  
\end{center}

\bigskip

\noindent 
F. F. Nichita\\
Institute of Mathematics "Simion Stoilow" of the Romanian Academy \\ 
P.O. Box 1-764, RO-014700 Bucharest, Romania \\
email: \textit{Florin.Nichita@imar.ro}\\

\noindent 
Bogdan Popovici\\
Horia Hulubei National Institute for Physics and Nuclear Engineering,\\
P.O.Box MG-6, Bucharest-Magurele, R0mania\\
email:\textit{ popobog@theory.nipne.ro}
\end{document}